\documentclass[a4paper,twoside,twocolumn]{article}

\usepackage[twocolumn,textwidth=18.31cm,columnsep=.5cm]{geometry}
\usepackage[utf8]{inputenc}
\usepackage{amsmath}
\usepackage{amsthm}
\usepackage{amssymb}
\usepackage[english]{babel}
\usepackage[numbers]{natbib}
\makeatletter

\newtheorem{thm}{Theorem}
\newtheorem{cor}{Corollary}
\newtheorem{defn}{Definition}

\begin{document}
\title{Flatness of Discrete-time Systems, a simple Approach}
\author{Schlacher Kurt\thanks{Institute of Automatic Control and Control Systems Technology, Johannes
Kepler University Linz, Altenberger Strasse 69, 4040 Linz, Austria (e-mail: \{kurt.schlacher, martin.lindorfer\}@jku.at)},  Lindorfer Martin\footnotemark[1]}
\date{}
\maketitle

\begin{abstract}
Flatness of discrete-time systems can be characterized by two simple
properties. There exists a map, a submersion, from the flat coordinates
and their forward shifts to the state and the input of the discrete-time
system, such that the system equations are fulfilled identically.
Flat coordinates, together with their shifts, also describe simple
shift systems. Therefore, this map transfers dynamic systems to dynamic
systems. Based on these facts, necessary conditions for a system,
to be flat, are derived. It is worth mentioning that methods from
differential geometry are not required. But these methods are used
to derive algorithms for the test, whether a system is flat or not,
and they are used to derive a flat parametrization.
\end{abstract}

\section{Introduction}

Flatness for lumped parameter systems has been introduced about 30
years ago, see e.g. \cite{Fliess}, \cite{Fliess2}. It became very
popular in the control community and has many applications in the
field of continuous-time systems. The problem of input to state linearization
by static feedback has been solved 40 years ago, see \cite{Jaku1},
but there is no equivalent final result for dynamic feedback known
until today.

Flatness of discrete-time systems can be defined analogously to the
continuous-time case, but one has the alternatives of forward shifts
or backward shifts to replace the time derivatives. The forward shift
is the commonly accepted choice. Therefore, it is used in this contribution,
too. Flatness for this class of systems is not so popular because
of the lack of applications. Nevertheless, it is an interesting problem
in system theory. See e.g. \cite{Jaku}, \cite{Bricaire_2}, \cite{Bricaire_1},
\cite{Kotta}, \cite{Kolar16_} for the development in this field.
One find newer approaches in \cite{Schlacher} or \cite{Kolar21},
\cite{Diwold_1}. The references of the last two ones give overviews
of the newer development in this field, too.

This contribution is an enhancement of the ideas presented in \cite{Schlacher}.
The discrete-time systems are identified with maps between manifolds.
Therefore, we recall some facts from differential geometry concerning
calculus on abstract manifolds in Section \ref{sec:NaP}. Section
\ref{sec:NF} presents a normal form for flat systems. Necessary and
locally sufficient conditions are given there for the existence of
a transformation to this form. In addition, it is shown, how one constructs
this transformation. The main result is presented in Section \ref{sec:F}.
The existence of the transformation from Section \ref{sec:NF} is
a necessary requirement for a system to be flat. It is worth mentioning
that properties of certain maps are exploited here, only. In contrast
to other contributions, differential geometry-based methods are not
used. The results of Section \ref{sec:F} and the conditions and methods
from Section \ref{sec:NF} lead to two approaches for tests to check,
whether a system is flat or not. These tests are presented in Section
\ref{sec:E}. The simpler one can be implemented in a computer algebra
system in a straightforward manner. If the system is flat, the result
is a system of PDEs for the flat outputs. The advanced test delivers
the flat parametrization if it exists. But it requires the solution
of linear PDEs or nonlinear ODEs.

\section{Notation and Prerequisites}

\label{sec:NaP}

Dynamic systems are modeled here by help of abstract manifolds. We
summarize the notation and some results concerning the differential
geometric methods used later on, for more information see e.g. \cite{kushner_lychagin_rubtsov_2006},
\cite{frankel_2011}. Let $\mathcal{M}$ be an $m$-dimensional smooth
manifold with local coordinates $x=x^{1},\ldots,x^{m}$. The set of
smooth functions $\text{\ensuremath{\mathcal{M}\rightarrow\mathbb{R}}}$
is denoted by $C^{\infty}\left(\mathcal{M}\right)$\footnote{The considerations here are not restricted to the smooth scenario,
but this assumption simplifies the notation.}. The functions $x^{1},\ldots,x^{m}$ are called canonical basis functions\footnote{The symbols $x^{i}$ are used for coordinates, as well as for the
functions $x^{i}$, where no confusion exists.}. Let the sequence $b^{i}\left(x\right)\in C^{\infty}\left(\mathcal{M}\right)$
be functionally independent, then $b^{1},\ldots,b^{m}$ is another
choice of basis functions. A bundle $\pi:\mathcal{S}\rightarrow\text{\ensuremath{\mathcal{X}}}$
is a triple with the total manifold $\mathcal{S}$, the base manifold
$\mathcal{X}$, and the surjective projection $\pi$. Locally, we
use the coordinates $x=x^{1},\ldots,x^{n}$ for $\mathcal{X}$ and
$x,u$, $u=u^{1},\ldots,u^{m}$ for $\mathcal{S}$ and get $\pi:x,u\mapsto x$.
The set of smooth sections or maps of the type $x=x$, $u=\sigma\left(x\right)$
is denoted by $\Gamma\left(\mathcal{S}\right)$. The tangent, cotangent
bundle of a manifold $\mathcal{M}$ are denoted by $\mathcal{T}\left(\mathcal{M}\right)$,
$\mathcal{T}^{\ast}\left(\mathcal{M}\right)$. Their standard bases
are $\left\{ \partial_{b^{1}},\ldots,\partial_{b^{m}}\right\} $,
$\left\{ \mathrm{d}b^{1},\ldots,\mathrm{d}b^{m}\right\} $ with $\left\langle \mathrm{d}b^{j},\partial_{b^{i}}\right\rangle =\delta_{i}^{j}$,
where $\delta_{i}^{j}$ denotes the Kronecker symbol. The exterior
derivative $\mathrm{d}:\Gamma\left(\bigwedge^{p}\left(\mathcal{M}\right)\right)\rightarrow\Gamma\left(\bigwedge^{p+1}\left(\mathcal{M}\right)\right)$
with $\bigwedge^{0}\left(\mathcal{M}\right)\cong C^{\infty}\left(\mathcal{M}\right)$,
$\bigwedge^{1}\left(\mathcal{M}\right)\rightarrow\mathcal{T}^{\ast}\left(\mathcal{M}\right)$
maps $p$-forms to $\left(p+1\right)$-forms, whereas the interior
product $\cdot\rfloor\cdot:\Gamma\left(\mathcal{T}\left(\mathcal{M}\right)\right)\times\Gamma\left(\bigwedge^{p+1}\left(\mathcal{M}\right)\right)\rightarrow\Gamma\left(\bigwedge^{p}\left(\mathcal{M}\right)\right)$
maps $\left(p+1\right)$-forms to $p$-forms.

A distribution $D=\mathrm{span}\left(B_{D}\right)$, $B_{D}=\left\{ v_{1},\ldots,v_{k}\right\} $,
$v_{i}\in\Gamma\left(\mathcal{T}\left(\mathcal{M}\right)\right)$
is a subspace $D_{x}\subset\mathcal{T}_{x}\left(\mathcal{M}\right)$,
or submodule $D\subset\mathcal{T}\left(\mathcal{M}\right)$, respectively.
It is called regular, iff its dimension is constant on an open neighborhood
of $x$. Iff $k$ is minimal, then $B_{D}$ is a basis of $D$. A
(regular) codistribution $D^{*}$ is defined analogously. We confine
the considerations to the regular case. The Lie bracket $\Gamma\left(\mathcal{T}\left(\mathcal{M}\right)\right)\times\Gamma\left(\mathcal{T}\left(\mathcal{M}\right)\right)\rightarrow\Gamma\left(\mathcal{T}\left(\mathcal{M}\right)\right)$
is denoted by $\left[\cdot,\cdot\right]$. A distribution is said
to be involutive, iff $\left[D,D\right]\subset D$ is met. The annihilator
$D^{\perp}$ of $D$ is defined by $\left\langle D^{\perp},D\right\rangle =\left\{ 0\right\} $.
Iff $D$ is involutive, then $D^{\perp}$ has a basis $\left\{ \mathrm{d}f^{k+1},\ldots,\mathrm{d}f^{m}\right\} $
of locally exact differentials. Let $D\subset E$ be two involutive
distributions, then one can show that there exists an involutive distribution,
a complement $D_{c}$, such that $D\oplus D_{c}=E$ is met. Let $P$
be a codistribution, the largest distribution $C$, which meets $C\rfloor P=\mathrm{span}\left(\left\{ 0\right\} \right)$
, $C\rfloor\mathrm{d}P\subset P$ is called the Cauchy characteristic
distribution. One can show that $C$ is involutive.

Let $f:\mathcal{M}\rightarrow\mathcal{N}$, $f=f^{1},\ldots,f^{n}$
be a smooth map between the manifolds $\mathcal{M}$, $\mathcal{N}$
with coordinates $x=x^{1},\ldots,x^{m}$, $y=y^{1},\ldots,y^{n}$.
The push forward of a tangent vector field $v\in\Gamma\left(\mathcal{T}\left(\mathcal{M}\right)\right)$
at $x\in\mathcal{M}$ is denoted by $f_{\ast}\left(v\right)\left(x\right)$.
Let $f$ be a diffeomorphism, then one gets the tangent vector field
$f_{\ast}\left(v\right)\circ f^{-1}\left(y\right)\in\Gamma\left(\mathcal{T}\left(\mathcal{N}\right)\right)$.
Let $f$ be a smooth submersion and $D$ a distribution on $\mathcal{N}$.
There exists a distribution $E$ on $\mathcal{M}$, such that $f_{\ast}^{-1}\left(D_{y}\right)=E_{x}$
for $y=f\left(x\right)$, see \cite{kushner_lychagin_rubtsov_2006}.
In the case $D=\mathcal{T}\left(\mathcal{N}\right)$, we choose functions
$f^{i}\left(x\right)$, $i=n+1,\ldots,m$, such that $f^{1}\left(x\right),\ldots,f^{m}(x)$
are basis functions of $\mathcal{M}$. Obviously, $f_{\ast}\left(v_{i}\right)=\partial_{y^{i}}$
with $v_{i}$ in $\left\{ \partial_{f^{1}},\ldots,\partial_{f^{n}}\right\} $
is met. The distribution $E=\mathrm{span}\left(B\right)$ with the
ring $f^{*}\left(C^{\infty}\left(\mathcal{N}\right)\right)$\footnote{The pull back of $p$-forms is denoted by $f^{\ast}:\bigwedge^{p}\left(\mathcal{N}\right)\rightarrow\bigwedge^{p}\left(\mathcal{M}\right)$.}
is a submodule of $\mathrm{span}\left(B\right)$. Elements of $E\oplus\mathrm{span}\left(\ker\left(\partial_{x}f\right)\right)$
are projectable tangent vector fields. They meet $f_{\ast}\left(v\right)\circ f^{-1}\left(y\right)\in\Gamma\left(\mathcal{T}\left(\mathcal{N}\right)\right)$.

\section{A Normal Form for Flat Systems}

\label{sec:NF}

A discrete-time system is a map $\mathcal{S}\rightarrow\mathcal{X}_{+}$,
where $\mathcal{S}$ is a bundle with the $n$-dimensional base manifold
$\mathcal{X}$, called state manifold, and the $\left(n+m\right)$-dimensional
total manifold $\mathcal{S}=\left(\mathcal{X},\mathcal{U}\right)$,
together with the projection $\pi:\mathcal{S}\rightarrow\mathcal{X}$.
The forward shift operator is denoted by $+$. We use the local coordinates
$x=x^{1},\ldots,x^{n}$ for $\mathcal{X}$ and $x,u$, $u=u^{1},\ldots,u^{m}$
for $\mathcal{S}$. In these coordinates, the system is represented
by 
\begin{equation}
x_{+}=f\left(x,u\right)=f^{1}\left(x,u\right),\ldots,f^{n}\left(x,u\right)\;,\label{eq:NF-01-a}
\end{equation}
with $f^{i}\in C^{\infty}\left(\mathcal{S}\right)$. $\mathcal{U}$
is called the input manifold, too. Sections of $\mathcal{S}$ are
control laws of the type $u=u\left(x\right)$. The bundle preserving
transformations are given by diffeomorphisms $t:\mathcal{S}\rightarrow\bar{\mathcal{S}}$
of the form 
\begin{equation}
\bar{x}=t_{x}\left(x\right)\;,\qquad\bar{u}=t_{u}\left(x,u\right)\;.\label{eq:NF-01-b}
\end{equation}
We assume that the functions $f^{{i}}$ are functionally independent.
Otherwise, the system is locally not reachable. In addition, the map
$f$ has constant rank on an open neighborhood of a point $\text{\ensuremath{\left(x,u\right)}}$,
or equivalently, $f$ is a submersion.

Let the map $f$ meet $\mathrm{rank}\left(\partial_{u}f\left(x,u\right)\right)=m_{u}$
on an open neighborhood of $\left(x,u\right)$ with $1\leq m_{u}\leq m$.
Possibly, after resorting the equations and variables, one can always
find functions $g^{i}=g\left(x,u\right)\in C^{\infty}\left(\mathcal{S}\right)$,
$i=m_{u}+1,\ldots,m$, such that the map 
\begin{equation}
\bar{x}=x\;,\qquad\bar{u}=t_{u}\left(x,u\right)\,,\label{eq:NF-02-a}
\end{equation}
see (\ref{eq:NF-01-b}), with $t^{i}\left(x,u\right)=f^{n-m_{u}+i}\left(x,u\right),$
$i=1,\ldots,m_{u}$ is invertible. The system takes the form 
\begin{equation}
\bar{x}_{+}=\bar{f}^{1}\left(\bar{x},\bar{u}\right),\ldots,\bar{f}^{n-m_{u}}\left(\bar{x},\bar{u}\right),\bar{u}^{1},\ldots,\bar{u}^{m_{u}}.\label{eq:NF-02-b}
\end{equation}
Obviously, $\bar{f}\left(\bar{x},\bar{u}\right)=\bar{f}^{1}\left(\bar{x},\bar{u}\right),\ldots,\bar{f}^{n-m_{u}}\left(\bar{x},\bar{u}\right)$
is independent of the trivial inputs $\bar{u}^{m_{u}+1},\ldots,\bar{u}^{m}$.
One can assign any value to these variables without impact on the
remaining system. Please note that $\mathrm{rank}$$\left(\partial_{x}f\left(x,u\right)\right)=n-m_{u}$
is met. To simplify the notation, we set $u=\bar{u}^{1},\ldots,\bar{u}^{m_{u}}$,
$m=m_{u}$ for the following considerations.

It is of interest to see, whether one can split the system in simpler
ones. We consider the transformation 
\begin{equation}
\hat{x}=x\;,\qquad\tilde{v},\tilde{w}=\hat{t}_{\tilde{v}}\left(x,u\right),\hat{t}_{\tilde{w}}\left(x,u\right)\label{eq:NF-03-a}
\end{equation}
with $\tilde{v}=\tilde{v}^{1},\ldots,\tilde{v}^{m_{v}}$, $\tilde{w}=\tilde{w}^{1},\ldots,\tilde{w}^{m_{w}},$
and $m_{v}+m_{w}=m$. The system takes the form 
\[
\hat{x}_{+}=\hat{f}\left(\hat{x},\tilde{v},\tilde{w}\right),\tilde{v},\tilde{w}\;,
\]
in the new coordinates. Roughly speaking, we replace $\mathcal{S}$
from above by its refinement $\mathcal{S}=\left(\mathcal{X},\mathcal{V}\times\mathcal{W}\right)$,
$\mathcal{U}=\mathcal{V}\times\mathcal{W}$. Let us assume, there
exists a submersion $\mathcal{S}\rightarrow\mathcal{X}_{+}$ 
\[
\tilde{x}_{+}=g\left(\hat{x},\tilde{v}\right),\tilde{v},\tilde{w}
\]
with $g=g^{1}\left(\hat{x},\tilde{v}\right),\ldots,g^{n-m}\left(\hat{x},\tilde{v}\right)$,
$g^{i}\in C^{\infty}\left(\left(\mathcal{X},\mathcal{V}\right)\right)$
and a diffeomorphism $h:\mathcal{X}_{+}\rightarrow\mathcal{X}_{+}$,
such that 
\begin{equation}
\hat{f}\left(\hat{x},\tilde{v},\tilde{w}\right),\tilde{v},\tilde{w}=h\left(g\left(\hat{x},\tilde{v}\right),\tilde{v},\tilde{w}\right)\label{eq:NF-03-b}
\end{equation}
is met. The state transformations $\hat{x}=h\left(\tilde{x}\right)$,
$\hat{x}_{+}=h\left(\tilde{x}{}_{+}\right)$ lead to 
\begin{equation}
\tilde{x}_{+}=g\left(\hat{x},\tilde{v}\right),\tilde{v},\tilde{w}=\tilde{g}\left(\tilde{x},\tilde{v}\right),\tilde{v},\tilde{w}=\tilde{f}\left(\tilde{x},\tilde{v},\tilde{w}\right)\;,\label{eq:NF-03-c}
\end{equation}
with $\tilde{g}\left(\tilde{x},\tilde{v}\right)=g\left(h\left(\tilde{x}\right),\tilde{v}\right)$.
Obviously, the system splits in the pure shift system 
\begin{equation}
\tilde{x}_{+}^{n-m_{w}+1},\ldots,\tilde{x}_{+}^{n}=\tilde{w}\label{eq:NF-04-a}
\end{equation}
and the remaining one 
\begin{equation}
\tilde{x}_{+}^{1},\ldots,\tilde{x}_{+}^{n-m_{w}}=g\left(\hat{x},\tilde{v}\right),\tilde{v}\label{eq:NF-04-b}
\end{equation}
with input $\tilde{v},\tilde{x}^{n-m_{w}+1},\ldots,\tilde{x}^{n}$.
The system (\ref{eq:NF-04-a}) describes a map $\mathcal{W}\rightarrow\mathcal{X_{+}}$.
The forward shift maps the functions $C^{\infty}\left(\mathcal{X}\right)$
to $C^{\infty}\left(\mathcal{X_{+}}\right)$ and is invertible in
this case. Combining these facts, the additional input $\tilde{x}^{n-m_{w}+1},\ldots,\tilde{x}^{n}$
of (\ref{eq:NF-04-a}) is easily derivable. In addition, these maps
can be extended to the tangent spaces $\mathcal{T}\left(\mathcal{W}\right),\mathcal{T}\left(\mathcal{X}\right),\mathcal{T}\left(\mathcal{X_{+}}\right)$
in a straight forward manner.

The special case $m_{v}=0$ is worth mentioning, since one gets 
\[
\tilde{x}_{+}=\tilde{g}\left(\hat{x}\right),w=\tilde{g}\left(h\left(\tilde{x}\right)\right),w\;.
\]
This is always the case for single input systems with $m=1$.

The system has particular properties, especially infinitesimal ones.
Let us introduce the system 
\begin{equation}
\tilde{P}_{0}=\mathrm{span}\left(\mathrm{d}g\left(\hat{x},\tilde{v}\right),\mathrm{d}\tilde{v},\mathrm{d}\tilde{w}\right)\label{eq:NF-05-a}
\end{equation}
together with the distributions 
\begin{align*}
\tilde{V} & =  \mathrm{span}\left(B_{\tilde{V}}\right)\;,\quad &B_{\tilde{V}}&=\left\{ \partial_{\tilde{v}^{1}},\ldots,\partial_{\tilde{v}^{m_{v}}}\right\} \\
\tilde{W} & =  \mathrm{span}\left(B_{\tilde{W}}\right)\;,\quad &B_{\tilde{W}}&=\left\{ \partial_{\tilde{w}^{1}},\ldots,\partial_{\tilde{w}^{m_{w}}}\right\} \;.
\end{align*}
We read off the following properties in a straightforward manner: 
\begin{enumerate}
\item $\tilde{f}_{*}\left(\mathrm{\partial}_{\tilde{w}^{i}}\right)=\partial_{\tilde{x}_{+}^{n+m_{v}+i}}\in\mathcal{T}\left(\mathcal{X}_{+}\right)$. 
\item Let $\tilde{P}_{1}\subset\tilde{P_{0}}$ be the largest codistribution
with $\tilde{W}\rfloor\tilde{P_{1}}=\mathrm{span}\left(\left\{ 0\right\} \right)$.
$\tilde{P}_{1}$ has a basis of exact differentials and $\tilde{W}\rfloor\mathrm{d}\tilde{P_{1}}=\mathrm{span}\left(\left\{ 0\right\} \right)$
is met. 
\item Iff $\tilde{f}_{*}\left(\tilde{v}\right)\in\mathcal{T}\left(\mathcal{X}_{+}\right)$,
$\tilde{v}\in\tilde{V}$ implies $\tilde{v}=0$, then $m_{W}$ is
maximal. 
\end{enumerate}
Since the systems (\ref{eq:NF-01-a}), (\ref{eq:NF-03-c}) are connected
by (\ref{eq:NF-01-b}), it remains to transfer these properties from
(\ref{eq:NF-05-a}) to the system 
\begin{equation}
P_{0}=\mathrm{span}\left(\mathrm{d}f\left(x,u\right)\right)\label{eq:NF-06-a}
\end{equation}
with $U=\mathrm{span}\left(B_{U}\right)$, $B_{U}=\left\{ \partial_{u^{1}},\ldots,\partial_{u^{m}}\right\} $,
such that they are necessary and locally sufficient for the existence
of a suitable transformation. 
\begin{enumerate}
\item There exists a distribution $W=\mathrm{span}\left(B_{W}\right)\subset U$,
$\dim$$\left(W\right)\geq1$ with a basis $B_{W}=\left\{ w_{1},\ldots,w_{m_{w}}\right\} $
and $f_{*}\left(w_{i}\right)\circ f^{-1}\left(x_{+}\right)\in\mathcal{T}\left(\mathcal{X}_{+}\right)$,
see Sec. \ref{sec:NaP}. 
\item Let $P_{1}\subset P_{0}$ be the largest codistribution with 
\begin{equation}
W\rfloor P_{1}=\text{span}(\left\{ 0\right\} )\,.\label{eq: W_floor_P}
\end{equation}
$P_{1}$ is integrable and $W$ is the Cauchy characteristic distribution
of $P_{1}$, which implies $W$ is involutive, see Sec. \ref{sec:NaP}. 
\item Since $U$ is involutive, there exists an involutive complement $V$
with $V\oplus W=U$, see Sec. \ref{sec:NaP}. Iff $v\not\in W$ is
projectable, implies $v=0$, then $m_{w}$ is maximal. 
\end{enumerate}
It remains to construct the transformation from (\ref{eq:NF-01-a})
to (\ref{eq:NF-03-c}). Since the system is integrable, there exist
functions $g^{i}\left(x,u\right)$, $i=1,\ldots,n-m_{w}$, such that
$P_{1}=\mathrm{span}\left(\mathrm{d}g^{1}\left(x,u\right),\ldots,\mathrm{d}g^{n-m_{w}}\left(x,u\right)\right)$
is met. Let us assume that the functions $g^{i}$, $f^{j}$ are sorted
such that the map $(u,v)=q\left(x,u\right)$ with 
\begin{eqnarray*}
v & = & g^{n-m}\left(x,u\right),\ldots,g^{n-m_{w}}\left(x,u\right)\\
w & = & f^{n-m_{w}+1}\left(x,u\right),\ldots,f^{n}\left(x,u\right)
\end{eqnarray*}
is invertible with respect to $u$. This is always possible because
of the construction of $V$, $W$. $W$ is the characteristic distribution
of $P_{1}$, which implies $\hat{g}^{i}\left(x,v\right)=g^{i}\left(x,q^{-1}\left(v,w\right)\right)$,
$i=1,\ldots,n-m$ or 
\[
P_{1}=\mathrm{span}\left(\mathrm{d}\hat{g}^{1}\left(x,v\right),\ldots,\mathrm{d}\hat{g}^{n-m}\left(x,v\right),\mathrm{d}v^{1},\ldots,\mathrm{d}v^{m_{v}}\right)
\]
is met. There must exist further functions $\hat{g}^{n-m+i}\left(x,v\right)$,
$i=1,\ldots,m$, such that the map $z=p\left(x,v\right)$ 
\[
z^{1},\ldots,z^{n}=\hat{g}^{1}\left(x,v\right),\ldots,\hat{g}^{n}\left(x,v\right)
\]
is invertible, otherwise the system (\ref{eq:NF-01-a}) is locally
not reachable. One can determine the map $h$, see (\ref{eq:NF-03-b}),
from the relation 
\[
f\left(p^{-1}\left(z,v\right),q^{-1}\left(v,w\right)\right)=h\left(z^{1},\ldots,z^{n-m},v,w\right)\;.
\]
Finally, the transformations $x=h\left(\tilde{x}\right)$, $x_{+}=h\left(\tilde{x}_{+}\right)$
lead to 
\[
\tilde{x}_{+}=h^{-1}\left(x_{+}\right)=\tilde{g}^{1}\left(\tilde{x},v\right),\ldots,\tilde{g}^{n-m}\left(\tilde{x},v\right),v,w
\]
with $\tilde{g}^{i}\left(\tilde{x},v\right)=\hat{g}^{i}\left(h\left(\tilde{x}\right),v\right)$,
$i=1,\ldots,n-m$.

\section{Flatness}

\label{sec:F}

Before we continue, we introduce the abbreviations $y_{a_{i},b_{i}}^{i}=y_{a_{i}}^{i},\ldots,y_{b_{i}}^{i}$,
$y_{a_{i},b_{i}}^{c,d}=y_{a_{c},b_{c}}^{c},\ldots,y_{a_{d},b_{d}}^{d}$
for sequences or sequences of sequences. 
\begin{defn}
The system (\ref{eq:NF-01-a}) is said to be locally flat (with respect
to forward shifts) iff there exists a smooth submersion $F=\left(F_{x},F_{u}\right)$
\begin{equation}
x=F_{x}\left(y_{0,r_{i}}^{1,m}\right)\;,\quad u=F_{u}\left(y_{0,r_{i}}^{1,m}\right)\;,\label{eq:F-01-a}
\end{equation}
such that 
\begin{equation}
F_{x,+}\left(y_{0,r_{i}}^{1,m}\right)=f\left(F_{x}\left(y_{0,r_{i}}^{1,m}\right),F_{u}\left(y_{0,r_{i}}^{1,m}\right)\right)\label{eq:F-01-b}
\end{equation}
is met. 
\end{defn}
The quantities $y_{0,r_{i}}^{1,m}$ are called flat variables with
their shifts. We assume that the number $r=\sum_{i=1}^{m}r_{i}$ is
minimal. Obviously, we can form $m$ simple shift systems of the type
\begin{equation}
y_{0,+}^{i}=y_{1}^{i},\ldots,y_{\left(r_{i}-1\right),+}^{i}=y_{r_{i}}^{i}\;,\quad i=1,\ldots,m\;.\label{eq:F-02-a}
\end{equation}
According to the assumptions of Section \ref{sec:NF}, $r_{i}\ge1$,
$i=1,\ldots,m$ is met.

The relation (\ref{eq:F-01-b}) can be rewritten as 
\begin{equation}
F_{x}\left(y_{1,r_{i}+1}^{1,m}\right)=f\left(F_{x}\left(y_{0,r_{i}}^{1,m}\right),F_{u}\left(y_{0,r_{i}}^{1,m}\right)\right)\;,\label{eq:F-01-c}
\end{equation}
or equivalently (\ref{eq:F-01-a}) is given by 
\begin{equation}
x=F_{x}\left(y_{0,r_{i}-1}^{1,m}\right)\;,\quad u=F_{u}\left(y_{0,r_{i}}^{1,m}\right)\;.\label{eq:F-01-d}
\end{equation}
If the system (\ref{eq:NF-01-a}) is transformed to (\ref{eq:NF-02-b}),
the relations (\ref{eq:F-01-d}) are converted to 
\begin{eqnarray*}
\bar{x} & = & \bar{F}_{\bar{x}}\left(y_{0,r_{i}-1}^{1,m}\right)=F_{x}\left(y_{0,r_{i}-1}^{1,m}\right)\\
\bar{u} & = & \bar{F}_{\bar{u}}\left(y_{1,r_{i}}^{1,m}\right)=t_{u}\left(F_{x}\left(y_{0,r_{i}-1}^{1,m}\right),F_{u}\left(y_{0,r_{i}}^{1,m}\right)\right)\;,
\end{eqnarray*}
see (\ref{eq:NF-02-a}). The relation (\ref{eq:F-01-c}) also shows
that the right-hand side must be independent of $y_{0}=y_{0}^{1},\ldots,y_{0}^{m}$.

The elimination of the variables $y_{r}=y_{r_{1}}^{1},\ldots,y_{r_{m}}^{m}$
from (\ref{eq:F-02-a}) leads to the subsystem 
\begin{equation}
y_{0,+}^{i}=y_{1}^{i},\ldots,y_{\left(r_{i}-2\right),+}^{i}=y_{r_{i}-1}^{i}\;,\quad i=1,\ldots,m\;.\label{eq:F-02-b}
\end{equation}
Since (\ref{eq:F-01-a}) maps the system (\ref{eq:F-02-a}) to (\ref{eq:NF-01-a}),
it must map (\ref{eq:F-02-b}) to a subsystem of (\ref{eq:NF-01-a}),
where some inputs $\tilde{w}$ are eliminated. This is possible only
if there exists a transformation 
\begin{eqnarray*}
\tilde{v} & = & \hat{F}_{\tilde{v}}\left(y_{1,r_{i}-1}^{1,m}\right)=\hat{t}_{\tilde{v}}\left(F_{x}\left(y_{0,r_{i}-1}^{1,m}\right),F_{u}\left(y_{0,r_{i}}^{1,m}\right)\right)\\
\tilde{w} & = & \hat{F}_{\tilde{w}}\left(y_{1,r_{i}}^{1,m}\right)=\hat{t}_{\tilde{w}}\left(F_{x}\left(y_{0,r_{i}-1}^{1,m}\right),F_{u}\left(y_{0,r_{i}}^{1,m}\right)\right)
\end{eqnarray*}
with $m_{w}=\mathrm{rank}\left(\partial_{y_{r}}F_{u}\left(y_{0,r_{i}}^{1,m}\right)\right)$,
see (\ref{eq:NF-03-a}) and (\ref{eq:F-01-d}). The assumption $r$
is minimal implies $m_{w}\geq1$. This subsystem has $n-m_{w}$ states
and $m_{v}+m_{w}$ inputs and must be of the type (\ref{eq:NF-03-c}).
Please note that some of the inputs can be trivial.

Let us consider the static subsystem $y_{0}=y_{0}^{1},\ldots,y_{0}^{m}$
of (\ref{eq:F-02-a}). One can assign any value to $y_{0}$ because
of lack of restricting conditions. Following the previous considerations,
there must exist an adequate $m$-dimensional subsystem of (\ref{eq:NF-01-a}).
This implies $\mathrm{rank}\left(\partial_{y_{0}}F_{x}\left(y_{0,r_{i}-1}^{1,m}\right)\right)=m$
must be met, otherwise, $r$ would not be minimal. Summarizing the
considerations of this section, we come to the following theorem. 
\begin{thm}
\label{thm:AnC}A necessary condition for the system (\ref{eq:NF-01-a})
to admit a flat parametrization (\ref{eq:F-01-a}) is, it must be
transformable to the form (\ref{eq:NF-03-c}). 
\end{thm}
It is worth mentioning that no tools from differential geometry are
used to derive this result. Only two properties have been used. The
map (\ref{eq:F-01-a}) must respect the system (\ref{eq:NF-01-a}),
or (\ref{eq:F-01-b}) must be fulfilled. In addition, (\ref{eq:F-01-a})
maps the dynamics of the pure shift system (\ref{eq:F-02-a}) to (\ref{eq:NF-01-a}).
This must apply to subsystems of (\ref{eq:F-02-a}), too.

Iff the system (\ref{eq:NF-01-a}) meets the condition of theorem
\ref{thm:AnC}, then one derives the flat subsystem (\ref{eq:NF-04-a}),
and the remaining one (\ref{eq:NF-04-b}). Now, one continues with
the remaining one until we end up with a trivial one or a non-flat
subsystem. The maximum number of repetitions is $n$. The following
result is a direct and simple consequence of this procedure since,
in each step, some of the inputs, see (\ref{eq:NF-04-b}), are replaced
by new ones. 
\begin{cor}
Let the system (\ref{eq:NF-01-a}) admit a flat parametrization (\ref{eq:F-01-a}).
Iff there are no trivial inputs, then the flat output is independent
of the input $u$.
\end{cor}

\section{Examples}

\label{sec:E}

We consider the example from \cite{Kolar16} given by 
\begin{equation}
\begin{aligned}x_{+}^{1} & =x^{2}(u^{1}+1)\\
x_{+}^{2} & =u^{1}\\
x_{+}^{3} & =x^{4}+u^{2}-1\\
x_{+}^{4} & =x^{5}+1-\frac{x^{1}(u^{1}+1)}{x^{2}+1}\\
x_{+}^{5} & =u^{2}+x^{2}.
\end{aligned}
\label{eq:system}
\end{equation}

\subsection{Transformations and Decompositions}

With the pfaffian system 
\begin{align*}
\begin{split}P_{0}=\, & \text{span}\left(\left\{ \vphantom{\frac{1}{x^{2}+1}}(u^{1}+1)\mathrm{d}x^{2}+x^{2}\mathrm{d}u^{1},\right.\right.\\
 & \left.\left.\mathrm{d}u^{1},\mathrm{d}x^{4}+\mathrm{d}u^{2},\mathrm{d}x^{2}+\mathrm{d}u^{2}\right.\right.,\\
 & \left.\left.\frac{1}{x^{2}+1}\left(-(u^{1}+1)\mathrm{d}x^{1}+\frac{(u^{1}+1)x^{1}}{x^{2}+1}\mathrm{d}x^{2}\right.\right.\right.\\
 & \left.\left.\left.\vphantom{\frac{1}{x^{2}+1}}\quad+(x^{2}+1)\mathrm{d}x^{5}-x^{1}\mathrm{d}u^{1}\right)\right\} \right)\,,
\end{split}
\end{align*}
see (\ref{eq:NF-06-a}), we derive the distributions $W_{0},V_{0}$
\begin{align*}
\begin{split}W_{0} & =\text{span}(\{\partial_{u^{2}}\})\\
V_{0} & =\text{span}(\{\partial_{u^{1}}\}).
\end{split}
\end{align*}
$P_{1}$, see Section \ref{sec:NF}, follows as 
\begin{align*}
\begin{split}P_{1}=\, & \text{span}\left(\left\{ \vphantom{\frac{1}{x^{2}+1}}\mathrm{d}u^{1},\mathrm{d}x^{2}-\mathrm{d}x^{4},(u^{1}+1)\mathrm{d}x^{2}+x^{2}\mathrm{d}u^{1}\right.\right.,\\
 & \left.\left.\frac{1}{x^{2}+1}\left(-(u^{1}+1)\mathrm{d}x^{1}+\frac{(u^{1}+1)x^{1}}{x^{2}+1}\mathrm{d}x^{2}\right.\right.\right.\\
 & \left.\left.\left.\vphantom{\frac{1}{x^{2}+1}}\quad+(x^{2}+1)\mathrm{d}x^{5}-x^{1}\mathrm{d}u^{1}\right)\right\} \right).
\end{split}
\end{align*}
Now, we are able to determine the functions $g(x,u)$ 
\begin{align*}
\begin{split}g^{1} & =\frac{-x^{1}\left(u^{1}+1\right)+x^{5}\left(x^{2}+1\right)}{x^{2}+1}\\
g^{2} & =x^{2}\left(u^{1}+1\right)\\
g^{3} & =x^{2}-x^{4}\\
g^{4} & =u^{1}.
\end{split}
\end{align*}
With the transformation $v^{1},w^{1}=q(x,u)=u^{1},u^{2}+x^{2}$, the
functions $\hat{g}(x,v)$ result in 
\begin{align*}
\begin{split}\hat{g}^{1} & =\frac{-x^{1}\left(v^{1}+1\right)+x^{5}\left(x^{2}+1\right)}{x^{2}+1}\\
\hat{g}^{2} & =x^{2}\left(v^{1}+1\right)\\
\hat{g}^{3} & =x^{2}-x^{4}.
\end{split}
\end{align*}
With these functions, one can determine the transformation ${x}=h(\tilde{x})$,
which looks as follows 
\begin{align*}
x^{1} & =\tilde{x}^{2},\quad & x^{2} & =\tilde{x}^{4},\quad & x^{3} & =\tilde{x}^{5}-\tilde{x}^{3}-1\\
x^{4} & =\tilde{x}^{1}+1,\quad & x^{5} & =\tilde{x}^{5},\quad &  & {}
\end{align*}
and the system in the new coordinates follows as 
\begin{align*}
\begin{split}\tilde{x}_{+}^{1} & =\frac{-\tilde{x}^{2}\left(v^{1}+1\right)+\tilde{x}^{5}\left(\tilde{x}^{4}+1\right)}{\tilde{x}^{4}+1}\\
\tilde{x}_{+}^{2} & =\tilde{x}^{4}\left(v^{1}+1\right)\\
\tilde{x}_{+}^{3} & =-\tilde{x}^{1}+\tilde{x}^{4}-1\\
\tilde{x}_{+}^{4} & =v^{1}\\
\tilde{x}_{+}^{5} & =w^{1}.
\end{split}
\end{align*}
Now, the last state can be eliminated. The states and the inputs will
be renamed as stated below 
\begin{align*}
\begin{split}x^{1} & =\tilde{x}^{1},\quad x^{2}=\tilde{x}^{2},\quad x^{3}=\tilde{x}^{3}\\
x^{4} & =\tilde{x}^{4},\quad u^{1}=\tilde{x}^{5},\quad u^{2}=v^{1}\,.
\end{split}
\end{align*}
This procedure will be repeated now and we derive the pfaffian system
\begin{align*}
\begin{split}P_{1}=\, & \text{span}\left(\left\{ \vphantom{\frac{1}{x^{2}+1}}\mathrm{d}u^{2},\mathrm{d}x^{4}-\mathrm{d}x^{1},(u^{2}+1)\mathrm{d}x^{4}+x^{4}\mathrm{d}u^{2}\right.\right.,\\
 & \frac{1}{x^{4}+1}\left((-u^{2}-1)\mathrm{d}x^{2}+\frac{(u^{2}+1)x^{2}}{x^{4}+1}\mathrm{d}x^{4}\right.\\
 & \left.\left.\left.\vphantom{\frac{1}{x^{4}+1}}\quad+(x^{4}+1)\mathrm{d}u^{1}-x^{2}\mathrm{d}u^{1}\right)\right\} \right)\,,
\end{split}
\end{align*}
the distributions 
\begin{align*}
\begin{split}W_{1} & =\text{span}(\{\partial_{u^{1}},\partial_{u^{2}}\})\\
V_{1} & =\text{span}(\{\})\,,
\end{split}
\end{align*}
and the reduced pfaffian system 
\begin{equation}
P_{2}=\text{span}(\{\mathrm{d}x^{4},\mathrm{d}x^{4}-\mathrm{d}x^{1}\})\,.\label{eq: P2}
\end{equation}
In this step, $m_{v}=0$ is met, and therefore the functions $g(x,u)$
result in 
\begin{align*}
\begin{split}g^{1} & =x^{4}\\
g^{2} & =x^{4}-x^{1}.
\end{split}
\end{align*}
Now, the functions $\hat{g}(x,v)$ remain unchanged and therefore
$\hat{g}(x,v)=g(x,u)$ holds. The input transformation is determined
by 
\begin{align*}
\begin{split}w^{1} & =\frac{({x}^{4}+1)u^{1}-({u}^{2}+1){x}^{2}}{x^{4}+1}\\
w^{2} & =u^{2}.
\end{split}
\end{align*}
The state transformation $x=h(\tilde{x})$ and the transformed system
result in 
\begin{align*}
x^{1} & =\tilde{x}^{3},\quad & x^{2} & =(\tilde{x}^{4}+1){\tilde{x}^{1}}\\
x^{3} & =\tilde{x}^{2}-1,\quad & x^{4} & =\tilde{x}^{4}
\end{align*}
and 
\begin{align*}
\tilde{x}_{+}^{1} & =\tilde{x}^{4},\quad & \tilde{x}_{+}^{2} & =\tilde{x}^{4}-\tilde{x}^{3}\\
\tilde{x}_{+}^{3} & =w^{1},\quad & \tilde{x}_{+}^{4} & =w^{2}\,.
\end{align*}
Since $m_{w}=2$ is met, two states can be eliminated and therefore,
the final system looks as follows 
\begin{align*}
\begin{split}x_{+}^{1} & =u^{2}\\
x_{+}^{2} & =u^{2}-u^{1},
\end{split}
\end{align*}
which shows that the system is flat.

\subsection{Transformations and Distributions}

Now, we check for flatness based on distributions. A similar result
can be found in \cite{Kolar19}, \cite{Kolar21}, derived by a different
mathematical machinery.

Section \ref{sec:NaP} shows, how one constructs the distribution
of all projectable vector fields. We choose the basis functions $f^{1},\ldots,f^{5}$
from (\ref{eq:system}) and $f^{6}=x^{1},f^{7}=x^{3}$. The vectors
\begin{align*}
\begin{split}\partial_{f^{1}} & =\frac{1}{u^{1}+1}\left(\partial_{x^{2}}+\partial_{x^{4}}-\partial_{u^{2}}\right)-\frac{x^{1}}{(x^{2}+1)^{2}}\partial_{x^{5}}\\
\partial_{f^{2}} & =\frac{x^{2}}{u^{1}+1}\left(\partial_{u^{2}}-\partial_{x^{2}}-\partial_{x^{4}}\right)-\frac{x^{1}(2x^{2}+1)}{(x^{2}+1)^{2}}\partial_{x^{5}}+\partial_{u^{1}}\\
\partial_{f^{3}} & =\partial_{x^{4}}\\
\partial_{f^{4}} & =\partial_{x^{5}}\\
\partial_{f^{5}} & =-\partial_{x^{4}}+\partial_{u^{2}}\\
\partial_{f^{6}} & =\partial_{x^{1}}+\frac{u^{1}+1}{x^{2}+1}\partial_{x^{5}}\\
\partial_{f^{7}} & =\partial_{x^{3}}
\end{split}
\end{align*}
form a basis $B=\left\{ \partial_{f^{1}},\ldots,\partial_{f^{7}}\right\} $
of $\mathcal{T}\left(\mathcal{S}\right)$, such that the first five
vectors are mapped to the unit vectors $\partial_{x_{+}^{1}},\ldots,\partial_{x_{+}^{5}}$.
A basis of the kernel of $f_{\ast}$ is given by 
\begin{align*}
\begin{split}B_{K}=\{\partial_{x^{3}},(x^{2}+1)\partial_{x^{1}}+(u^{1}+1)\partial_{x^{5}}\}\;.\end{split}
\end{align*}
The vector fields $\partial_{u^{1}},\partial_{u^{2}}$ are given in
the basis $B$ as follows 
\begin{align*}
\begin{split}\partial_{u^{1}} & =x^{2}\partial_{f^{1}}+\partial_{f^{2}}-\frac{x^{1}}{x^{2}+1}\partial_{f^{4}}\\
\partial_{u^{2}} & =\partial_{f^{3}}+\partial_{f^{5}}\,.
\end{split}
\end{align*}
Obviously $\partial_{u^{2}}$ is projectable. But no nontrivial combination
with $\partial_{u^{1}}$ is projectable because of $x^{1}/(x^{2}+1)\notin f^{\ast}\left(C^{\infty}\left(\mathcal{X}\right)\right)$.

A computationally simpler method is given by the following schema.
Given $B_{U}=\left\{ \partial_{u^{1}},\ldots,\partial_{u^{m}}\right\} $,
we set $f_{\ast}\left(U\right)=\mathrm{span}\left(\left\{ f_{\ast}\left(\partial_{u^{1}}\right),\ldots,f_{\ast}\left(\partial_{u^{m}}\right)\right\} \right)$
and construct $V=\left\{ v_{1}\left(x,u\right),\ldots,v_{m}\left(x,u\right)\right\} $
with normalized elements such that $f_{\ast}\left(U\right)=\mathrm{span}\left(V\right)$
is met. Normalization means, there are $m$ indices $j_{1},\ldots,j_{m}$
such that $v_{i}^{j_{k}}\left(x,u\right)=\delta_{i}^{k}$ is fulfilled.
A field $v=\sum_{i=1}^{m}\lambda^{i}v_{i}$ is projectable, iff $v^{i}\left(x,u\right)\in f^{\ast}\left(C^{\infty}\left(\mathcal{X}\right)\right)$
is met. Because of the normalization, we get $\lambda^{i}\in f^{\ast}\left(C^{\infty}\left(\mathcal{X}\right)\right)$.
Therefore, the equations 
\begin{equation}
\mathrm{d}v^{i}\wedge\mathrm{d}f^{1}\wedge\cdots\wedge\mathrm{d}f^{n}=0\label{eq: Form_check}
\end{equation}
are linear equations in $\lambda^{i}$, which allows us to determine
a basis of the projectable vector fields. Possibly, one has to repeat
the procedure with the set of solutions for $\lambda_{1},\ldots,\lambda_{m}$.

Now we apply the procedure form above to example (\ref{eq:system}).
The pushforward of the vector fields $\partial_{u^{1}}$, $\partial_{u^{2}}$
result in 
\[
\begin{aligned}f_{*}(\partial_{u^{1}}) & =x^{2}\partial_{x_{+}^{1}}+\partial_{x_{+}^{2}}-\frac{x^{1}}{x^{2}+1}\partial_{x_{+}^{4}}=v_{1}\\
f_{*}(\partial_{u^{2}}) & =\partial_{x_{+}^{3}}+\partial_{x_{+}^{5}}=v_{2}\,.
\end{aligned}
\]
We determine $\mathrm{vol}=\mathrm{d}f^{1}\wedge\dots\wedge\mathrm{d}f^{n}$,
and according to (\ref{eq: Form_check}), we get the equations 
\begin{alignat*}{3}
 & \lambda^{1}\underbrace{\mathrm{d}x^{2}\wedge\mathrm{vol}}_{=0} &  & =0\qquad\qquad\lambda^{1}\underbrace{\mathrm{d}1\wedge\mathrm{vol}}_{=0} &  & =0\\
 & \lambda^{2}\underbrace{\mathrm{d}1\wedge\mathrm{vol}}_{=0} &  & =0\qquad\frac{-\lambda^{1}}{x^{2}+1}\underbrace{\mathrm{d}x^{1}\wedge\mathrm{vol}}_{\neq0} &  & =0\\
 & \lambda^{2}\underbrace{\mathrm{d}1\wedge\mathrm{vol}}_{=0} &  & =0\,,
\end{alignat*}
with a solution 
\[
\lambda^{1}=0,\quad\quad\lambda^{2}=1\,.
\]
Now, we can determine the distributions 
\begin{align*}
W_{0}=\text{span}(\{\partial_{u^{2}}\}),\quad V_{0}=\text{span}(\{\partial_{u^{1}}\})\,.
\end{align*}
The codistribution $P_{1}$, see (\ref{eq: W_floor_P}), follows as
\begin{align*}
P_{1}=\, & \text{span}\left(\left\{ \vphantom{\frac{1}{x^{2}+1}}\mathrm{d}u^{1},\mathrm{d}x^{2}-\mathrm{d}x^{4},(u^{1}+1)\mathrm{d}x^{2}+x^{2}\mathrm{d}u^{1}\right.\right.,\\
 & \left.\left.\frac{1}{x^{2}+1}\left(-(u^{1}+1)\mathrm{d}x^{1}+\frac{(u^{1}+1)x^{1}}{x^{2}+1}\mathrm{d}x^{2}\right.\right.\right.\\
 & \left.\left.\left.\vphantom{\frac{1}{x^{2}+1}}\quad+(x^{2}+1)\mathrm{d}x^{5}-x^{1}\mathrm{d}u^{1}\right)\right\} \right).
\end{align*}
With the projectable vector field 
\[
f_{*}(\partial_{u^{2}})=\partial_{x_{+}^{5}}+\partial_{x_{+}^{3}}\,,
\]
we get the distribution 
\[
U_{1}=\text{span}(\{\partial_{x^{5}}+\partial_{x^{3}},\partial_{u^{1}}\}),
\]
see (\ref{eq:NF-04-a}), (\ref{eq:NF-04-b}). Repeating the procedure,
in the next steps, we derive $U$, $V$, $W$, and $P$ as 
\begin{align*}
\begin{split}W_{1} & =\text{span}(\{\partial_{x^{5}}+\partial_{x^{3}},\partial_{u^{1}}\})\\
V_{1} & =\text{span}(\{\})\\
P_{2} & =\text{span}(\{\mathrm{d}x^{2},\mathrm{d}x^{2}-\mathrm{d}x^{4}\})
\end{split}
\end{align*}
and 
\begin{align*}
\begin{split}U_{2} & =\text{span}\left(\left\{ \partial_{x^{4}},x^{1}\partial_{x^{1}}+(x^{2}+1)\partial_{x^{2}}\right\} \right)\\
W_{2} & =\text{span}\left(\left\{ \partial_{x^{4}},x^{1}\partial_{x^{1}}+(x^{2}+1)\partial_{x^{2}}\right\} \right)\\
V_{2} & =\text{span}(\{\})\,.
\end{split}
\end{align*}
$P_{2}$ is apart from state transformations identical to (\ref{eq: P2}).
Now, $\text{dim}(W_{0}\oplus W_{1}\oplus W_{2})=n$ is met, and the
distribution $\overline{W}=W_{0}\oplus W_{1}\oplus W_{2}$ results
in 
\begin{align*}
\overline{W}=\text{span}\left(\left\{ \partial_{u^{2}},\partial_{x^{5}}+\partial_{x^{3}},\partial_{u^{1}},\partial_{x^{4}},\right.\right.\\
\quad\left.\left.x^{1}\partial_{x^{1}}+(x^{2}+1)\partial_{x^{2}}\right\} \right)\,.
\end{align*}
The flat outputs are annihilated by $\overline{W}$, and follow as
\[
y^{1}=\frac{x^{2}+1}{x^{1}},\quad y^{2}=x^{5}-x^{3}.
\]

\section{Summary}

The main result of this contribution is a necessary condition for
discrete-time systems to be flat. The system must admit the transformation
to a normal form. This form allows a reduction to a flat system and
a complement. It is worth mentioning that only two simple facts are
used. A flat parametrization is a submersion from the flat coordinates
and their forward shifts to the state and the input of the discrete-time
system such that the parametrized system equations are fulfilled identically.
Since flat coordinates, together with their shifts, also describe
a transformation from the simple shift systems to the flat system,
this map must be valid for subsystems, too. It is worth mentioning
that tools from differential geometry are not necessary for the derivation
of this result.

The presented tests, whether a system is transformable to the normal
form, use methods from differential geometry. The simpler one checks
if a system is flat. It only requires the solution of linear equations
as well as substitution and simplification in nonlinear expressions.
If the system is flat, the result is a system of PDEs. Flat outputs
are a solution of these equations. It is worth comparing the simplicity
of this test with other ones. The advanced test produces a sequence
of systems such that a successor is a true subsystem of its predecessor.
The result of this test is the flat parametrization if it exists.
But it requires the solution of linear PDEs or nonlinear ODEs.

\bibliographystyle{plainnat}
\bibliography{linschla_ref}
\end{document}